\documentclass[reqno,b5paper]{amsart}
\usepackage{amsmath}
\usepackage{amssymb}
\usepackage{amsthm}
\usepackage{enumerate}
\usepackage[mathscr]{eucal}

\newcommand{\C}{\mathcal{C}}

\newcommand{\te}{{\tau c}}
\newcommand{\bh}{B({\mathcal{H}})}
\newcommand{\h}{{\mathcal{H}}}
\newtheorem{claim}{claim}[section]
\newtheorem{theorem}[claim]{Theorem}

\newtheorem{lemma}[claim]{Lemma}
\newtheorem{proposition}[claim]{Proposition}

\theoremstyle{definition}

\newtheorem{example}[claim]{Example}

\begin{document}
\title[Skew-symmetric operators and reflexivity]{Skew-symmetric operators and reflexivity}
\author[Chafiq Benhida \and Kamila Kli\'s-Garlicka \and Marek Ptak]{Chafiq Benhida* \and Kamila Kli\'s-Garlicka**\and Marek Ptak**,***}
\newcommand{\acr}{\newline\indent}
\address{\llap{*\,}Universite Lille 1\acr
        Laboratoire Paul Painlev\'{e} \acr
        59655 Villeneuve d'Ascq\acr France}
\email{chafiq.benhida@math.univ-lille1.fr}

\address{\llap{**\,}Department of Applied Mathematics\acr
University of Agriculture\acr Balicka 253c\acr 30-198 Krakow\acr Poland}\email{rmklis@cyf-kr.edu.pl}

\address{\llap{***\,}Institut of Mathematics\acr Pedagogical University\acr  ul. Podchor\c a\.zych 2\acr 30-084 Krak\'ow\acr Poland}
\email{rmptak@cyf-kr.edu.pl}

\thanks{The first named author was partially supported by Labex CEMPI (ANR-11-LABX-0007-01). The research of the second and the third author was financed by the Ministry of Science and Higher Education of the Republic of Poland}


\begin{abstract}
In contrast to the subspaces of all $C$-symmetric operators, we show that the subspaces of all skew-C symmetric operators are reflexive and even hyperreflexive with the  constant
$\kappa(\C^s)\leqslant 3$.
\end{abstract}

\maketitle

\section{Introduction and Preliminaries}
Let $\h$ be  a complex Hilbert space with an inner product   $\langle ., .\rangle$ and let
$\bh$ be the Banach algebra of all bounded linear operators on $\h$.

Recall that the space of trace class operators $\te$ is
predual to $\bh$ with the dual action $\langle T,
f\rangle=tr(Tf)$, for $T\in\bh$ and $f\in\te$. The trace norm in
$\te$ will be denoted by $\|\cdot\|_1$. By $F_k$ we denote the set of all
operators of rank at most $k$.  Often rank-one operators are written as $x\otimes y$,
for $x,\,y\in\h$, and $(x\otimes y)z=\langle z,y\rangle x$ for
$z\in\h$. Moreover, $tr(T(x\otimes y))=\langle Tx,y\rangle$.
Let $\mathcal{S}\subset B(\h)$ be a closed subspace. Denote by $\mathcal{S}_\bot$ the {\it preanihilator} of $\mathcal{S}$, i.e., $\mathcal{S}_\bot=\{t\in \tau c:  tr(St)=0 \text{ for all } S\in\mathcal{S}\}$.
A weak${}^*$ closed subspace $\mathcal{S}$ is {\it $k$-reflexive} iff rank-k operators  are linearly dense in $\mathcal{S}_\bot$, i.e., $\mathcal{S}_\bot=[\mathcal{S}_\bot\cap F_k]$ (see \cite{ls}). $k$-hyperreflexivity introduced in \cite{Ar1, KKMP} is a stronger property than $k$-reflexivity, i.e., each $k$-hyperreflexive subspace is $k$-reflexive. A subspace $\mathcal{S}$ is called {\it $k$-hyperreflexive} if there is a constant $c>0$ such that
\begin{equation}\label{hyp}\operatorname{dist} (T, \mathcal{S})\leq c\cdot \sup\{|tr(Tt)|: t\in F_k\cap \mathcal{S}_\bot, \|t\|_1\leq 1\},\end{equation}
for all $T\in B(\h)$. Note that $\operatorname{dist} (T, \mathcal{S})$ is the infimum distance. The supremum on the right hand side of \eqref{hyp} will be denoted by $\alpha_k(T,\mathcal{S})$. The smallest constant for which inequality \eqref{hyp} is satisfied is called the $k$-{\it hyperreflexivity constant} and is denoted $\kappa_k(\mathcal{S})$. If $k=1$, the letter $k$ will be omitted.

Recall that $C$ is a {\it conjugation} on $\h$ if $C:\h\longrightarrow\h$ is an antilinear,  isometric involution, i.e.,
$\langle Cx,Cy\rangle =\langle y,x\rangle$ for all $x, y\in \h$ and $C^2=I$.
   An operator $T$ in $\bh$ is said to be {\it $C$--symmetric} if  $CTC=T^*$.
   $C$--symmetric operators have been intensively studied by many authors in the last decade (see \cite{Ga 3}, \cite{Ga 4}, \cite{GW}, \cite{K-P}).
   It is a wide class of operators including Jordan blocs, truncated Toeplitz operators and Hankel operators.

   Recently, in \cite {K-P}, the authors  considered the problem of reflexivity and hyperreflexivity  of the subspace
   $\C=\{{T}\in \bh : CTC=T^*\}$. They have shown that $\C$ is transitive and 2-hyperreflexive.
Recall that $T\in \bh$ is  a {\it skew--$C$ symmetric} iff $CTC=-T^*$.
 In this paper, $\C^s=\{{T}\in \bh : CTC=-T^*\}$  --
    the subspace of all skew--$C$ symmetric operators  will be investigated from the reflexivity and hyperreflexivity point of view. It follows directly from the definition that $\C$ and $\C^s$ are weak${}^*$ closed.

  We emphasize that the notion of skew symmetry  is  linked to many problems in physics and that
  any operator $T \in \bh$ can be written as a sum of a $C$--symmetric operator and a skew--$C$ symmetric operator.
  Indeed, $T=A+B$, where  $A=\frac{1}{2}(T+CT^*C)$ and $B=\frac{1}{2}(T-CT^*C)$.

  The aim of this paper is to show that $\C^s$ is  reflexive and even hyperreflexive.

\section{Preanihilator}
Easy calculations show the following.
\begin{lemma}\label{1} Let $C$ be a conjugation in a complex Hilbert space $\h$ and  $h,g\in\h$. Then
\begin{enumerate}
  \item  $C(h \otimes g)C=Ch\otimes Cg$,
  \item $h \otimes g-Cg \otimes Ch \in \C^s$ .
\end{enumerate}
\end{lemma}

In \cite[Lemma 2]{Ga 3} it was shown that \[\C \cap F_1 = \{\alpha\cdot h \otimes Ch : h \in \h, \alpha\in \mathbb{C}\}.\]
We will show that it is also a description of the rank-one operators in the preanihilator of $\C^s$.
\begin{proposition}\label{2} Let $C$ be a conjugation in a complex Hilbert space $\h$. Then
\[\C^s_{\perp} \cap F_1 = \C \cap F_1= \{\alpha\cdot h \otimes Ch : h \in \h, \alpha\in \mathbb{C}\}. \]
\end{proposition}

\begin{proof}
To prove "$\supset$" let us take $T\in \C^s$ and $h \otimes Ch \in \C \cap F_1$. Then
\begin{multline*}
\langle T, h \otimes Ch \rangle = \langle Th, Ch\rangle  = \langle h, CTh\rangle = \langle h, -T^*Ch \rangle = \\ = - \langle Th, Ch\rangle= - \langle T, h \otimes Ch \rangle.
\end{multline*}
Hence $\langle T, h \otimes Ch \rangle=0$ and $ h \otimes Ch \in \C^s_{\perp} \cap F_1$.

For the converse inclusion let us take a rank-one operator $ h \otimes Cg \in \C^s_{\perp}$. Since $Cg \otimes h - Ch \otimes g \in \C^s$, by Lemma \ref{1} we have
\begin{multline*}
0=\langle Cg \otimes h -Ch \otimes g, h\otimes Cg\rangle = \langle (Cg \otimes h)h, Cg\rangle  - \langle (Ch \otimes g) h , Cg\rangle =\\=  \|h\|^2 \cdot \|Cg\|^2 - \langle h,g \rangle \langle Ch, Cg\rangle=  \|h\|^2 \cdot \|g\|^2 - |\langle h,g \rangle|^2.
\end{multline*}

Hence $|\langle h,g \rangle|= \|h\| \, \|g\|$, i.e., there is equality in Cauchy-Schwartz inequality. Thus $h, g$ are linearly dependent and the proof in finished.
\end{proof}

\begin{lemma}\label{3}
Let $C$ be a conjugation in a complex Hilbert space $\h$. Then
\[\C^s_{\perp} \cap F_2 \supset \{h \otimes g + Cg \otimes Ch : h,g \in \mathcal{H}\}. \]
\end{lemma}

\begin{example} Note that for different conjugations we obtain different subspaces.
 Let $C_1(x_1, x_2, x_3)=(\bar x_3,\bar x_2, \bar x_1)$ be a conjugation on $\mathbb{C}^3$. Then
 $$\C_1^s=\left\{\left(
    \begin{array}{ccc}
      a & b & 0 \\
      c & 0 & -b \\
      0 & -c & -a \\
    \end{array}
  \right): a, b, c\in\mathbb{C}\right\}
$$
and
$$\mathcal{C}_1=\left\{ \left(
                       \begin{array}{ccc}
                         a & b & \ast \\
                        c & \ast & b \\
                        \ast & c & a \\
                       \end{array}
                     \right): a, b, c\in\mathbb{C}\right\}.
$$
Rank-one operators in $\mathcal{C}_1$ and in $(\C^s_1)_\bot$ are of the form $\alpha(x_1, x_2, x_3)\otimes(\bar x_3,\bar x_2, \bar x_1)$ for $\alpha\in\mathbb{C}$.

 If we now consider another conjugation $C_2(x_1, x_2, x_3)=(\bar x_2,\bar x_1, \bar x_3)$ on $\mathbb{C}^3$, then
 $$\C_2^s=\left\{\left(
    \begin{array}{ccc}
      a & 0 & b \\
      0 & -a & c \\
     -c & -b & 0 \\
    \end{array}
  \right): a, b, c\in\mathbb{C}\right\},
 $$
and $$ \mathcal{C}_2=\left\{ \left(
                       \begin{array}{ccc}
                         a & \ast & b \\
                        \ast & a & c \\
                         c & b & \ast \\
                       \end{array}
                     \right): a, b, c\in\mathbb{C}\right\}.
$$ Similarly, rank-one operators in $\mathcal{C}_2$ and in $(\C^s_2)_\bot$ are of the form $\alpha(x_1, x_2, x_3)\otimes(\bar x_2,\bar x_1, \bar x_3)$.
\end{example}

\begin{example}
Let $C$ be a conjugation in $\h$. Consider
$\tilde{C}=\left(  \begin{array}{cc}0&C\\C&0\end{array}\right)$ the conjugation in $\h\oplus \h$ (see \cite{LZ}). An operator $T\in B(\h\oplus\h)$ is skew-$\tilde{C}$ symmetric, if and only if
$T=\left(  \begin{array}{cc} A& B\\ D& -CA^*C
\end{array}                                                                                                                                \right)$,
    where  $A, B, D\in B(\h)$ and $B$, $D$
are skew-C symmetric. Moreover, rank-one operators in $\mathcal{\tilde{C}}^s_\bot$   are of the form $\alpha (f\oplus g)\otimes(Cg\oplus Cf)$ for $f, g\in \mathcal{H}$ and $\alpha\in\mathbb{C}$.
\end{example}

\begin{example}
Let us consider the classical Hardy space $H^2$ and let $\alpha$ be a nonconstant inner function. Define $K^2_\alpha=H^2\ominus \alpha H^2$ and $C_\alpha h(z)=\alpha \overline{zh(z)}$. Then $C_\alpha$ is a conjugation on $K^2_\alpha$. By $S_\alpha$ and $S^*_\alpha$ denote the compressions of the unilateral shift $S$ and the backward shift $S^*$ to $K^2_\alpha$, respectively.
Recall after \cite{Sarason} that the kernel functions in $K^2_\alpha$ for $\lambda\in \mathbb{C}$ are projections of appropriate kernel functions $k_\lambda$ onto $K^2_\alpha$, namely $k^\alpha_\lambda=k_\lambda-\overline{\alpha(\lambda)}\alpha k_\lambda$. Denote by $\tilde{k}_\lambda^\alpha=C_\alpha k_\lambda^\alpha$.
Since $S_\alpha$ and $S^*_\alpha$ are $C_\alpha$ -- symmetric (see \cite{Ga 3}),  for a skew--$C_\alpha$ symmetric operator $A\in B(K^2_\alpha)$ we have
\begin{multline}\label{4}
\langle AS^n_\alpha k_\lambda^\alpha, (S^*_\alpha)^m \tilde{k}_\lambda^\alpha \rangle=\langle C_\alpha (S^*_\alpha)^m\tilde{k}_\lambda^\alpha, C_\alpha AS^n_\alpha k_\lambda^\alpha\rangle=\\-\langle S_\alpha^m C_\alpha\tilde{k}_\lambda^\alpha,A^*C_\alpha S^n_\alpha k_\lambda^\alpha  \rangle=- \langle A S_\alpha^m k_\lambda^\alpha, (S_\alpha^*)^n \tilde{k}_\lambda^\alpha\rangle,
\end{multline}
for all $n, m\in \mathbb{N}$. Note that  if $n=m$, then
\begin{equation}\label{toepl}\langle A S_\alpha^n k_\lambda^\alpha, (S_\alpha^*)^n \tilde{k}_\lambda^\alpha\rangle=0.\end{equation}

In particular, we may consider the special case $\alpha=z^k$, $k>1$. Then the equality \eqref{toepl} implies that a skew--$C_\alpha$ symmetric operator $A\in B(K^2_{z^k})$ has the matrix representation in the canonical basis with $0$ on the diagonal orthogonal to the main diagonal. Indeed, let $A\in B(K^2_{z^k})$ have the matrix $(a_{ij})_{i,j=0,\dots,k-1}$ with respect to the canonical basis. Note that $C_{z^k}f=z^{k-1}\bar f$, $k_0^{z^k}=1$, $\tilde{k}_0^{z^k}=z^{k-1}$. Hence for $0\leq n\leq k-1$ we have
$$0=\langle A S_\alpha^n 1, (S_\alpha^*)^n z^{k-1}\rangle=\langle Az^n,z^{k-n-1}\rangle=a_{n,k-n-1}.$$
Moreover, from the equality \eqref{4} we can obtain that
$$\langle Az^n,z^{k-m-1}\rangle=-\langle Az^m,z^{k-n-1}\rangle, $$ which implies that $a_{n,k-m-1}=-a_{m,k-n-1}$ for $0\leq m, n\leq k-1$.
\end{example}

\section{Reflexivity}
The following theorem can be obtained as a corollary of Theorem \ref{hyp}. However, we think that the proof presented here is also interesting.
\begin{theorem} Let $C$ be a conjugation in a complex Hilbert space $\h$.
The subspace $\C^s$ of all skew--C symmetric operators on $\h$ is reflexive.
\end{theorem}
\begin{proof}
By Proposition \ref{2}  it is necessary to show that if $\langle T, h \otimes Ch \rangle=\langle Th, Ch \rangle=0$ for any $h \in \h$, then $CTC=-T^*$.

Recall after \cite[Lemma 1]{Ga 3}  that $\mathcal{H}$ can be decomposed into its real and imaginary parts $\mathcal{H}=H_{R}+i\,H_I$.
  Recall also that we can write  $h=h_R+ih_I\in \mathcal{H}$ with $h_R=\tfrac 12(I+C)h \in H_R$ and $h_I=\tfrac 1{2i}(I-C)h \in H_I$. Then $Ch_R=h_R$, $Ch_I=h_I$ and $Ch=C(h_R+ih_I)=h_R-ih_I$.

Let $T \in \bh$. The operator $T$ can be represented as
$\left[\begin{array}{cc}
W & X \\
Y& Z \\\end{array}\right]$,
where $W \colon H_R \to H_R$, $Z \colon H_I \to H_I$, $X \colon H_I \to H_R$, $Y \colon H_R \to H_I$ and $W,\;X,\;Y,\;Z$ are real linear. The condition $CTC=-T^*$ is equivalent to the following: $W=-W^*$, $Z=-Z^*$, $Y=X^*$.

On the other hand, the condition $\langle Th, Ch \rangle=0$ for any $h=h_R+ih_I$ is equivalent to
\begin{equation}\label{11}
\langle Wh_R,h_R\rangle + \langle Xh_I, h_R\rangle - \langle Yh_R,h_I\rangle -\langle Zh_I, h_I \rangle=0
\end{equation}
for any $h_R \in H_R$, $h_I \in H_I$.
In particular, $\langle Wh_R,h_R\rangle=0$ for any $h_R \in H_R$.

Let $h^{'}_R, h^{''}_R\in H_R$. Then
$\langle Wh^{'}_R,h^{'}_R\rangle=0$, $\langle Wh^{''}_R,h^{''}_R\rangle=0$ and \[0=\langle W(h^{'}_R+h^{''}_R), h^{'}_R+h^{''}_R \rangle =
\langle Wh^{'}_R, h^{''}_R\rangle + \langle Wh^{''}_R, h^{'}_R \rangle.\]
Hence \[\langle Wh^{'}_R, h^{''}_R\rangle=\langle h^{'}_R, -Wh^{''}_R\rangle\]
and finally $W^*=-W$.
Since, by (\ref{11}), in particular $\langle Zh_I,h_I \rangle=0$ for any $h_I\in H_I$ we can also get $Z^*=-Z$.

Because $\langle Wh_R,h_R \rangle=0=\langle Zh_I,h_I \rangle$ for any $h_R\in H_R$, $h_I\in H_I$,  hence by \eqref{11} we get \[\langle Xh_I,h_R \rangle-\langle Yh_R,h_I \rangle=0.\] Thus $Y=X^*$ and the proof is finished.
\end{proof}
Recall that a single operator $T\in B(\h)$ is called reflexive if the weakly closed algebra generated by $T$ and the identity is reflexive.  In  \cite{LZ} authors characterized normal skew symmetric operators and by  \cite{Sar} we know that every normal operator is reflexive. Hence one may wonder, if all skew--C symmetric operators are reflexive. The following simple example shows that it is not true.
\begin{example}
Consider the space $\mathbb{C}^2$ and a conjugation $C(x,y)=(\bar x, \bar y)$. Note that operator
$T=\left(
     \begin{array}{cc}
       0 & 1 \\
       -1 & 0 \\
     \end{array}
   \right)
$ is skew-C symmetric. The weakly closed algebra $\mathcal{A}(T)$ generated by $T$ consists of operators of the form $\left(
                                                                                                       \begin{array}{cc}
                                                                                                        a & b \\
                                                                                                         -b & a \\
                                                                                                       \end{array}
                                                                                                     \right)$. Hence $\mathcal{A}(T)_\bot=\left\{ \left(
                                                                                                                                               \begin{array}{cc}
                                                                                                                                                 t & s \\
                                                                                                                                                 s & -t \\  \end{array}\right): t,s\in\mathbb{C} \right\}$. It is easy to see, that $\mathcal{A}(T)_\bot\cap F_1=\{0\}$, which implies that $T$ is not reflexive.
\end{example}

\section{Hyperreflexivity}
\begin{theorem}\label{hyp1} Let $C$ be a conjugation in a complex Hilbert space $\h$. Then the subspace
$\C^s$   of all skew--$C$ symmetric operators is hyperreflexive with the constant $\kappa(\C^s)\leqslant 3$ and $2$-hyperreflexive with $\kappa_2(\C^s)=1$.
\end{theorem}
\begin{proof}
Let $A \in \bh$. Firstly, similarly as in the proof of Theorem 4.2 \cite{K-P} it can be shown that $A-CA^*C\in \C^s$. It is also shown there that $(CAC)^*=CA^*C$. Hence we have
\begin{align*}
&d(A, \C^s)\leqslant \|A-\tfrac 12(A-CA^*C)\|=\tfrac 12\|A+CA^*C\|
\\&=\tfrac 12 \|A^*+CAC\|= \tfrac 12 \sup\{|\langle h, (A^*+CAC)g\rangle| : \|h\| , \|g\| \leqslant 1\}
\\&= \tfrac 12 \sup\{|\langle A, h \otimes g + Cg \otimes Ch\rangle| : \|h\| , \|g\| \leqslant 1\}
\\&= \tfrac 12 \sup\{|\langle A, h \otimes Cg + g \otimes Ch\rangle| : \|h\|, \|g\| \leqslant 1\}=\alpha_2(A, \C^s)
\\&= \tfrac 12 \sup\{|\langle A, (h+g) \otimes C(h+g)-h \otimes Ch - g \otimes Cg\rangle| : \|h\|, \|g\| \leqslant 1\}
\\ &\leqslant \tfrac 12 \sup\{|\langle A, h \otimes Ch\rangle| : \|h\| \leqslant 1\} +  \tfrac 12 \sup\{|\langle A, g \otimes Cg\rangle| : \|g\| \leqslant 1\} +
 \\& +  \tfrac 12 \sup\{4|\langle A, \tfrac 12(h+g)  \otimes C(\tfrac 12(h+g))\rangle| : \|h\|, \|g\| \leqslant 1\}\\
&\leqslant 3\  \alpha (A,\C^s).
\end{align*}
We have used the characterization given in Proposition \ref{2}.
\end{proof}


\begin{thebibliography}{99}

\bibitem{Ar1}
W. T. Arveson:
\textit {Interpolation problems in nest algebras}, J. Funct. Anal.  \textbf{ 20}  (1975), 208--233.
\bibitem{Ga 3} S. R. Garcia and M. Putinar: \textit{Complex symmetric operators and applications}, Trans. Amer. Math. Soc.  \textbf{358} (2006), 1285--1315.
\bibitem{Ga 4} \underline{\hspace{2cm}}: \textit{Complex symmetric operators and applications} II, Trans. Amer. Math. Soc.  \textbf{359} (2007), 3913--3931.
\bibitem{GW} S. R. Garcia and W. R. Wogen: \textit{Some new classes of complex symmetric operators}, Trans. Amer. Math. Soc.  \textbf{362} (2010), 6065--6077.
\bibitem{K-P} K. Kli\'s-Garlicka and M. Ptak: \textit{C-symmetric operators and reflexivity}, Operators and Matrices \textbf{9} no. 1 (2015),  225--232.
\bibitem{KKMP} K. Kli\'s, M. Ptak: \textit{$k$-hyperreflexive subspaces}, Houston J. Math.   \textbf{32} (2006), 299--313.
\bibitem{LZ} C. G. Li and S. Zhu: \textit{Skew symmetric normal operators}, Proc. Amer. Math. Soc.  \textbf{141} no. 8 (2013), 2755--2762.
\bibitem{ls} A.I. Loginov and V.S. Shul'man: \textit{Hereditary and
intermediate reflexivity of W*-algebras}, Izv. Akad. Nauk.
SSSR,  \textbf{39} (1975), 1260--1273; Math. USSR-Izv. \textbf{9}  (1975), 1189--1201
\bibitem{Sarason} D. Sarason: \textit{Algebraic properties of truncated Toeplitz operators}, Operators and Matrices \textbf{1}  no. 4 (2007), 491--526.
\bibitem{Sar} \underline{\hspace{1,5cm}}: \textit{Invariant subspaces and unstarred operator algebras}, Pacific J. Math.  \textbf{17} no. 3 (1966), 511--517.
\bibitem{Zu} S. Zhu: \textit{Skew symmetric weighted shifts}, Banach J. Math. Anal.  \textbf{9} no. 1 (2015), 253--272.
\bibitem{Zu1} \underline{\hspace{1cm}}:  \textit{Approximate unitary equivalence to skew symmetric operators}, Complex Anal. Oper. Theory  \textbf{8} no. 7 (2014), 1565--1580.
\end{thebibliography}
\end{document}